\documentclass[12pt]{article}
\usepackage{amsfonts}
\begin{document}

\def\href#1#2{#2}
\newpage
\pagestyle{empty}
\setcounter{page}{0}

\begin{center}
{\huge LAPTH} \\
\bigskip
%%{\large \em Laboratoire d'Annecy-Le-Vieux de Physique Th\'eorique}
\end{center}
\smallskip
\hrule

\vspace{20mm}

\begin{center}

  {\LARGE  {\sffamily 
      A Note on the Generalised Lie Algebra $sl(2)_q$
      }}\\[1cm]
  
  \vspace{10mm}
  
  {\large Daniel Arnaudon \footnote{Daniel.Arnaudon@lapp.in2p3.fr.}}\\[.42cm]
  
  {\em Laboratoire  d'Annecy-Le-Vieux de Physique
    Th{\'e}orique }LAPTH\footnote{URA 1436 
    du CNRS, associ{\'e}e 
    {\`a} l'Universit{\'e} de Savoie.
    \newline\indent
    Work partially supported by European Community contract TMR
    FMRX-CT96.0012}, CNRS\\[.242cm]
  
  LAPP, BP 110, F-74941
  Annecy-le-Vieux Cedex, France.
  \\

\end{center}
\vspace{20mm}

\begin{abstract}
  In a recent paper, V. Dobrev and A. Sudbery classified the
  highest-weight and lowest-weight finite dimensional irreducible 
  representations
  of the quantum Lie algebra $sl(2)_q$ introduced by V. Lyubashenko and
  A. Sudbery. The aim of this note is to add to this classification
  all the finite dimensional irreducible representations which have no
  highest weight and/or no lowest weight, in the case when $q$ is a
  root of unity. 
  For this purpose, we give a description of the enlarged centre.
\end{abstract}

\rightline{math.QA/9804131}
\rightline{LAPTH-679/98}
\rightline{April 1998}

\newpage
\pagestyle{plain}
%%%%%%%%%%%%%%%%%%%%%%%%%%%%%%%%%%%%%%%%%%%%%%%%%%%%%%%%%%%%%%%%%%%%

%\def\mathbb{\bf}
\def\CC{{\mathbb C}}
\def\NN{{\mathbb N}}
\def\QQ{{\mathbb Q}}
\def\RR{{\mathbb R}}
\def\ZZ{{\mathbb Z}}
\def\cA{{\cal A}}          \def\cB{{\cal B}}          \def\cC{{\cal C}}
\def\cD{{\cal D}}          \def\cE{{\cal E}}          \def\cF{{\cal F}}
\def\cG{{\cal G}}          \def\cH{{\cal H}}          \def\cI{{\cal I}}
\def\cJ{{\cal J}}          \def\cK{{\cal K}}          \def\cL{{\cal L}} 
\def\cM{{\cal M}}          \def\cN{{\cal N}}          \def\cO{{\cal O}}
\def\cP{{\cal P}}          \def\cQ{{\cal Q}}          \def\cR{{\cal R}} 
\def\cS{{\cal S}}          \def\cT{{\cal T}}          \def\cU{{\cal U}}
\def\cV{{\cal V}}          \def\cW{{\cal W}}          \def\cX{{\cal X}}
\def\cY{{\cal Y}}          \def\cZ{{\cal Z}}
%
%%\bibliographystyle{utphys}
%%\bibliographystyle{unsrt}
%\def\href#1#2{#2}
%
% D{\'e}finitions locales:
%
\def\uq#1{{\cal U}_q(sl(#1))}
\renewcommand{\thefootnote}{\fnsymbol{footnote}}

\section{\label{sect:introduction}Introduction}

In the notion of ``quantum groups'' introduced by Drinfeld and Jimbo
\cite{Dri,Jim}, one actually refers to the quantisation of the enveloping
algebra $\cU(\cG)$, considered as a Hopf algebra. 
The question arises then about the existence of a deformation of the
Lie algebra itself, and several authors have more recently made progresses 
towards a definition of quantised  Lie algebras
\cite{Majid:jgp94,DeliusHuff:jpa29,LyuSud:1995}.   

In \cite{DobSud}, V.~Dobrev and A.~Sudbery give 
a classification of finite dimensional irreducible
representations of the quantum Lie algebra $sl(2)_q$ as defined in 
\cite{LyuSud:1995} by V.~Lyubashenko and A.~Sudbery. 
This classification actually concerns the highest
weight and lowest weight representations. It happens however that there
exists other classes of finite dimensional representations 
of quantum groups at roots of unity that are
useful for physics, namely the periodic (cyclic) representations
\cite{Skl,RA}, which 
appear for instance in
generalisations of the chiral Potts model \cite{BazKas:cmp136,DJMMa}.

The quantum Lie algebra is defined 
as a finite dimensional subspace of the quantised enveloping
algebra that is invariant under the quantum adjoint action. 
According to \cite{LyuSud:1995}, the representation theory of
$sl(2)_q$ reduces to that of the albegras $\cB$ and $\cF$ defined
below. 

The algebra $\cB$ is generated by $X_0$, $X_\pm$, $C$, related by 
\begin{eqnarray}
  &&  q^2 X_0 X_+ - X_+ X_0 = q C X_+ \;,\\
  &&  q^{-2} X_0 X_- - X_+ X_0 = - q^{-1} C X_- \;,\\
  &&  X_+ X_- - X_- X_+ = (q+q^{-1}) (C-\lambda X_0) X_0 \;,\\
  &&  C X_\pm - X_\pm C = C X_0 - X_0 C = 0 \;,
  \label{eq:relB}
\end{eqnarray}
where $\lambda=q-q^{-1}$. We will later use $q$-numbers $[p]$ defined
as usual 
by $[p]\equiv\frac{\displaystyle q^p-q^{-p}}{\displaystyle q-q^{-1}}$.

A quadratic central element of $\cB$ is given by 
\begin{equation}
  C'_2 = X_- X_+ + q^{-1} C X_0 + q^{-2} X_0^2
  \label{eq:C2p}
\end{equation}
(normalised by a factor $(q+q^{-1})^{-1}$ with respect to \cite{DobSud}.)

The algebras $\cF$ and $\cA$ are defined from $\cB$ by adding
respectively the relations 
$C^2 - \lambda^2 C'_2 = 1$ and $C=1$
on central elements \cite{LyuSud:1995}.

When interpreted in the $\uq2$ context, $C$ corresponds to the usual
quadratic Casimir element, whereas the quadratic central element
$C'_2$ of $\cB$ corresponds to a quartic central element
\cite{LyuSud:1995}. 

\bigskip
The classification of finite dimensional irreducible representations
of $\uq2$ at roots of unity (including periodic ones) was given in
\cite{RA}. The classification we present here is very close to the
latter. The representations of $\cB$ have one more parameter. Unusual
representations of dimension 1 are present.

The methods we use here are similar to those used 
in \cite{ABosp12q} in the case of
$\cU_q(osp(1|2))$, and more details may be found there.

\section{\label{sect:centre}Centre at $q^{2l}=1$}

A linear basis of $\cB$  is given by 
\begin{equation}
  X_-^{a_-} X_+^{a_+} X_0^{a_0} C^{b_1} {C'_2}^{b_2}
  \qquad \mbox{with} \qquad a_\pm, a_0, b_1, b_2 \in \NN, \qquad 
  a_+ a_- = 0 \;.
  \label{eq:basis}
\end{equation}
This can be proved starting from a basis of the form 
given in \cite{LyuSud:1995}, Lemma 3.2.  
Then all the common powers of $X_-$ and $X_+$ in a monomial 
can indeed be reexpressed 
in terms involving $X_0$, $C$ and $C'_2$ only, using 
\begin{eqnarray}
  X_- X_+ &=& C'_2 - q^{-1} C X_0 - q^{-2} X_0^2 
  \label{eq:XmXp}
  \\
  &=& \lambda^{-2} 
  \Big\{ - (C^2 - \lambda^2 C'_2) + (1+q^{-2})C(C-\lambda X_0)
  - q^{-2} (C-\lambda X_0)^2 \Big\} \;, 
  \label{eq:XmXp2}
  \\
  X_+ X_- &=& C'_2 + q C X_0 - q^{2} X_0^2 
  \label{eq:XpXm}
  \\
  &=& \lambda^{-2} 
  \Big\{ - (C^2 - \lambda^2 C'_2) + (1+q^{2})C(C-\lambda X_0)
  - q^{2} (C-\lambda X_0)^2 \Big\}
  \;.
  \label{eq:XpXm2}
\end{eqnarray}

The centre of $\cB$ for generic $q$ is generated by $C$ and $C'_2$.
(A linear combination of terms given by (\ref{eq:basis}) needs, in
order to commute with $X_0$, to involve only terms with $a_+=a_-=0$. In
order to commute with $X_\pm$, it should not involve terms with
$a_0\neq 0$.)

\bigskip

\paragraph{Let now $q$ be a root of unity.} 
More precisely let $l$ be the smallest
(non 0) integer such that $q^{2l}=1$.

The centre of $\cB$ is now $\CC[C,C'_2,X_+^l,(C-\lambda
X_0)^l]+\CC[C,C'_2,X_-^l,(C-\lambda X_0)^l]$. The sum is not a direct
sum and the intersection is $\CC[C,C'_2,(C-\lambda X_0)^l]$. 

The generators $C$, $C'_2$, $X_\pm^l$ and
$(C-\lambda X_0)^l$ of the centre of $\cB$ when $q^{2l}=1$ are 
subject to the relation 
\begin{equation}
  X_-^l X_+^l = q^{l(l-1)} \lambda^{-2l}
  \left\{     
%%    - ( C^2 - \lambda^{2} C'_2 )^l
    - ( \cD^2 )^l
    +  q^{-l} \cD^l \cQ_l((q+q^{-1})C\cD^{-1}) (C-\lambda X_0)^{l} 
    - (C-\lambda X_0)^{2l} 
  \right\}\;,
  \label{eq:rel-centre}
\end{equation}
where $\cQ_l$ is the polynomial of degree $l$, related to the
Chebichev polynomial of the first kind, such that 
\begin{equation}
  \cQ_l(x+x^{-1}) = x^l + x^{-l} \;,
  \label{eq:Chebychev}
\end{equation}
and where $\cD$ is defined by 
\begin{equation}
  \cD^2 = C^2 - \lambda^2 C'_2\;. 
  \label{eq:Dcarre}
\end{equation}
Note
that the right hand side of (\ref{eq:rel-centre}) is a well defined
polynomial of degree $l$ in $\cD^2$, and hence in $C'_2$.

To prove the formula (\ref{eq:rel-centre}), we proceed as in
\cite{KerDipl}: we first prove by a simple
recursion 
\begin{equation}
  X_-^p X_+^p = \lambda^{-2p} \prod_{r=0}^{p-1} 
  q^{-2r-1}
  \Big\{
  - q\cD^2 q^{2r}
  + (q+q^{-1}) C (C-\lambda X_0) 
  - q^{-1} (C-\lambda X_0)^2 q^{-2r}
  \Big\} \;,
  \label{eq:rel-recurs}
\end{equation}
and then let $p=l$, so that the operand runs over all the powers of
$q^2$.

\section{Finite dimensional irreducible representations of $\cB$}

We now give the classification of finite dimensional irreducible
representations when $q^{2l}=1$, insisting on those with no highest
weight (and/or lowest weight) vector, which were not considered in
\cite{DobSud}. We use module notations.

On any finite dimensional simple module, the central
elements $C$, $C'_2$, $X_\pm^l$ and
$(C-\lambda X_0)^l$ act as scalars (diagonal matrices with a single
eigenvalue), which we denote respectively 
by $c$, $c'_2$, $x_\pm^l$ and $z$, and which satisfy the relation
(obtained from (\ref{eq:rel-centre}))
\begin{equation}
  x_-^l x_+^l = q^{l(l-1)} \lambda^{-2l}
  \left\{     
%%    - ( c^2 - \lambda^{2} c'_2 )^l
    - ( d^2 )^l
    +  q^{-l} d^l \cQ_l((q+q^{-1})c d^{-1}) z
    - z^{2} 
  \right\}\;,
  \label{eq:rel-centre-scal}
\end{equation}
where $d^2 \equiv c^2 - \lambda^2 c'_2$. Note that
(\ref{eq:rel-centre-scal}) is a polynomial of degree $l$ 
in $d^2$, and hence in $c'_2$. 

Let $M$ be a finite dimensional simple module. There exists in $M$
a vector $v_0$ such that, in addition to 
$C v_0 = c v_0$, 
$C'_2 v_0 = c'_2 v_0$,
$X_\pm^l v_0 = x_\pm^l v_0$ and 
$(C-\lambda X_0)^l v_0 = z v_0$, we also have 
\begin{itemize}
\item $X_0 v_0 = x_0 v_0$ with $z = (c-\lambda x_0)^l$,
\item $M=\mbox{span}\{X_+^p v_0, X_-^p v_0\}_{p=0,...,l-1}$ (these
  vectors being linearly dependent).
\end{itemize}
The existence of $v_0$ satisfying the first property is guaranteed by
the finite dimension. The second property is proved by writing 
$M=\cB .v_0$, using the basis (\ref{eq:basis}), and observing that
$X_\pm^p v_0$ are eigenvectors of $C$, $C'_2$ and $X_0$.

\paragraph{First case: $z\neq 0$ and $x_- \neq 0$ ($X_-$ acts
  injectively).} 

We define 
\begin{equation}
  v_p = x_-^{-p} X_-^p v_0 \qquad (v_l\equiv v_0)\;.
  \label{eq:vp}
\end{equation}
Then 
\begin{eqnarray}
  X_0 v_p &=& \left(q^{2p} x_0 - q^{p}[p] c \right) v_p \;,
  \label{eq:X0vp}
  \\
  X_- v_p &=& x_- v_{p+1} \;.
  \label{eq:X-vp}
\end{eqnarray}
\begin{equation}
    X_+ v_p = x_-^{-1} \lambda^{-2}
    \Big\{
    - d^2 %%=(c^2 - \lambda^2 c'_2)
    + (1+q^{-2})c(c-\lambda x_0) q^{2p}
    - q^{-2} (c-\lambda x_0)^2 q^{4p} 
    \Big\}  v_{p-1}\;.
  \label{eq:X+vp}
\end{equation}
The action of $X_+$ on $v_p$ is computed using 
$  X_+ v_p = x_-^{-1} X_+ X_- v_{p-1} $
and Eq. (\ref{eq:XpXm2}).
The module spanned by $v_p$, $p=0,...,l-1$ 
is simple since the eigenvectors $v_p$ of $X_0$ 
correspond to $l$ different eigenvalues 
(this would not be the case with $z=0$).

This class of periodic (or semi-periodic when $x_+=0$) $l$-dimensional 
representations is then  characterised by the five complex
parameters $c$, $c'_2$, $x_\pm^l$ and $z$, related by the polynomial
relation (\ref{eq:rel-centre-scal}).

\paragraph{Second case: $z \neq 0$, $x_- = 0$ and $x_+ \neq 0$ ($X_+$
  acts injectively, but not $X_-$).}
This case is symmetric to a subcase of the previous one, for which
$x_+=0$ was not excluded. 

Let $w_0=v_0$ (such that $X_0 w_0 = x_0 w_0$) 
requiring further that  $X_-w_0=0$. Such a vector exists
because i) $X_-$ is nilpotent ii) the eigenspace of $X_-$ related to
the eigenvalue 0 is stable under the action of $X_0$. 
We have $c'_2=q^2x_0^2-qcx_0$ and $z=(c-\lambda x_0)^l$.
We define 
\begin{equation}
  w_p = x_+^{-p} X_+^p w_0 \qquad (w_l\equiv w_0)\;.
  \label{eq:wp}
\end{equation}
Then
\begin{eqnarray}
  X_0 w_p &=& \left(q^{-2p} x_0 + q^{-p}[p] c \right) w_p \;,
  \label{eq:X0wp}
  \\
  X_+ w_p &=& x_+ w_{p+1} \;.
  \label{eq:X+wp}
\end{eqnarray}
\begin{equation}
    X_- w_p = x_+^{-1} \lambda^{-2}
    \Big\{
    - d^2 %%=(c^2 - \lambda^2 c'_2)
    + (1+q^{2})c(c-\lambda x_0) q^{-2p}
    - q^{2} (c-\lambda x_0)^2 q^{-4p} 
    \Big\}  w_{p-1} \;.
  \label{eq:X-wp}
\end{equation}
The action of $X_-$ on $w_p$ is computed using 
$  X_- w_p = x_+^{-1} X_- X_+ w_{p-1} $
and Eq. (\ref{eq:XmXp2}).
The module spanned by $w_p$, $p=0,...,l-1$ is again simple.

This class of semi-periodic $l$-dimensional 
representations is then  characterised by the three complex
parameters $c$, $x_0$, $x_+^l$. 
The parameters $c'_2$ and $z$ are related to those by 
$c'_2=q^2x_0^2-qcx_0$ and $z=(c-\lambda x_0)^l$.

\paragraph{Third case: $z\neq 0$ and  $x_-=x_+=0$ (Highest weight and
  lowest weight representation).} 
This case has been treaded in details in \cite{DobSud}. We just give a
summary of the classification given there. 

\begin{itemize}
\item 
There are one parameter irreducible representations of dimension
$n<l$, described by  
\begin{eqnarray}
  C v_p &=& c v_p \;,\\
  X_0 v_p &=& \lambda^{-1} \left(c - q^{2p} \nu \right) v_p \;,\\
  X_- v_p &=& v_{p+1} \;, \qquad\qquad X_- v_{n-1}=0 \;,\\
%%  X_+ v_p &=& q^{n+2p-2} \nu^2 [p][n-p] v_{p-1}
  X_+ v_p &=& \lambda^{-1} [p]q^{p-2} \nu 
  \left\{ (q^2+1)c - (q^{2p}+1) \nu  \right\} v_{p-1} \;,
  \label{eq:dim=n}
\end{eqnarray}
with the constraint $(q^2+1)c = (q^{2n}+1) \nu$
\footnote{Note that with this parametrisation, it is not necessary to
  distinguish the case $q^{2n}+1=0$, i.e. $n=l/2$, when $l$ is even,
  for which $c=0$. This is however not true for representations of
  $\cA$, for which $c=1$.} 
and $\nu\neq 0$. Note
that when $l=2$, $n=1$, this is a two parameter representation.

\item
There are $l$-dimensional irreducible representations, 
also described by (\ref{eq:dim=n}) (with $n=l$), and 
characterised by two parameters $c$ and $\nu\neq 0$, with the
constraint that $(q^2+1)c - (q^{2p}+1) \nu \neq 0$ for
$p=1,...,l-1$. These representations do not exist when $l=2$. 
\end{itemize}

\paragraph{Fourth case: $z=0$.}
Supposing first $x_-\neq 0$, we define $v_p$, $p=0,...,l-1$ as in the
first case. The action of $X_0$, $X_\pm$ are as in
(\ref{eq:X0vp},\ref{eq:X-vp},\ref{eq:X+vp}). Now, this defines a reducible
representation since all the eigenvalues of $X_0$ are
equal. Irreducible one dimensional subrepresentations are defined by
any vector $v=\sum_{p=0}^{l-1}q^{2kp}v_p$, and
\begin{equation}
  \begin{array}{l}
    C v = c v ,\cr
    C'_2 v = c'_2 v ,\cr
    X_0 v = x_0 v ,\cr
    X_\pm v = x'_\pm v ,\cr
  \end{array}  \qquad \mbox{with} \qquad  
  \begin{array}{l}
    c-\lambda x_0=0 ,\cr
    x'_\pm = q^{\pm 2k} x_\pm ,\cr
    x_+ x_- =  c'_2 - \lambda^{-2} c^2 = - \lambda^{-2} d^2 .
  \end{array}
  \label{eq:dim1}
\end{equation}
Considering then the case $x_-=0$, $x_+\neq 0$, 
and following the construction defined by 
(\ref{eq:wp},\ref{eq:X0wp},\ref{eq:X+wp},\ref{eq:X-wp}) again leads to 
(\ref{eq:dim1}). The case $z=x_-=x_+=0$, already in the classification of
\cite{DobSud}, is also described by (\ref{eq:dim1}). 

This class of one-dimensional representations described by 
(\ref{eq:dim1}) is characterised by three continuous 
parameters $x_0$, $x'_\pm$. 

\bigskip

\paragraph {Remark:} Even in the case when $q$ is generic, there exists 
(semi)-periodic representations of dimension 1, given by
\begin{equation}
  X_0 v = x_0 v\;, \qquad C v = c v \;, \qquad X_\pm v = x_\pm v \;,
  \qquad \mbox{with} \qquad  c-\lambda x_0=0 \;.
  \label{eq:dim1generique}
\end{equation}

\paragraph{Representations of $\cA$: } 
The irreducible finite dimensional representations of $\cA$ are given
by fixing $c=1$ in the previous classification. This is in general
possible, except for the representations of dimension $l/2$, (when
$l/2\in\NN$)  for which the constraint was $c=0$.

\section{Finite dimensional irreducible representations of $\cF$}

The algebra $\cF$ is defined from $\cB$ as its quotient by the
relation $C^2 = 1+\lambda^2 C'_2$, i.e. $\cD^2=1$ (\ref{eq:Dcarre}).
One obtains the irreducible finite dimensional representations of
$\cF$ from those of $\cB$ by imposing the supplementary condition
$d^2=c^2-\lambda^2 c'_2=1$ on the parameters. 
Generically, the parameters are then $c$, $x_\pm^l$ and $z$,
eigenvalues of $C$, $X_\pm^l$ and $(C-\lambda X_0)^l$, related by 
\begin{equation}
  x_-^l x_+^l = q^{l(l-1)} \lambda^{-2l}
  \left\{     
%%    - ( c^2 - \lambda^{2} c'_2 )^l
    - 1
    +  q^{-l} d^l \cQ_l((q+q^{-1})c) z
    - z^{2} 
  \right\} \;.
  \label{eq:rel-centre-scal-F}
\end{equation}
We still consider only the case when $q$ is a root of unity. 

\medskip
\noindent
The classification is then the following:
\paragraph{First case: $z\neq 0$ and $x_- \neq 0$.} The
representations with injective action of $X_-$, of dimension $l$,
described by  (\ref{eq:vp},\ref{eq:X0vp},\ref{eq:X-vp},\ref{eq:X+vp})
with $d^2=1$.  
They depend on the parameters $c$,
$x_\pm^l$ and $z$, related by (\ref{eq:rel-centre-scal-F}).

\paragraph{Second case: $z \neq 0$, $x_- = 0$ and $x_+ \neq 0$.}
The representations with nilpotent action of $X_-$ and injective
action of $X_+$, of dimension $l$, described by
(\ref{eq:wp},\ref{eq:X0wp},\ref{eq:X+wp},\ref{eq:X-wp}) with $d^2=1$. 
This class of semi-periodic $l$-dimensional 
representations  depends on the parameters $\nu$, 
$x_+^l$, from which $c$, $x_0$ and $z$ are given by 
$c=(q\nu+q^{-1}\nu^{-1})/[2]$,
$x_0=\lambda^{-1}(c-\nu)$ and $z=\nu^l$. 

\paragraph{Third case: $z\neq 0$ and  $x_-=x_+=0$.} 
This case corresponds to the classification in \cite{DobSud}. 
\begin{itemize}
\item 
The representations of dimension $n<l$ are described by
(\ref{eq:dim=n}) with $[2]c=q^{-1}\nu+q\nu^{-1}$ and
$\nu^2=q^{-2n+2}$. Hence, they are labeled by the dimension 
$n$ and a sign
$\epsilon$ such that $\nu=\epsilon q^{-n+1}$. 

\item
The representations of dimension $l$ are described by (\ref{eq:dim=n})
with again $[2]c=q^{-1}\nu+q\nu^{-1}$, and now $\nu^2\neq q^{-2p+2}$
for $p=1,...,l-1$. They are labeled by one parameter $\nu$. 
\end{itemize}

\paragraph{Fourth case: $z=0$.}
The unusual representations of dimension 1 described by
(\ref{eq:dim1}) still exist for $\cF$, with $d^2=1$. 
These representations are necessarily
periodic since $x_+ x_- = -\lambda^{-2}$, which explains that they are
not in the classification of \cite{DobSud}. 
They depend on two parameters $x_0$ and $x_+$.

\bigskip

\begingroup\raggedright\endgroup

\end{document}